\documentclass{article}
\usepackage{amscd}
\begin{document}

\pagestyle{myheadings}

\markright{MARCELO GOMEZ MORTEO}

\title{A Dirac Morphism for the Farrell-Jones Isomorphism Conjecture in K-Theory }
\author{ Marcelo Gomez Morteo }

\maketitle
\vspace{16pt}

\begin{abstract}

\vspace{16pt}

We construct a Dirac Morphism. We prove that if this Dirac morphism is invertible, then the isomorphism conjecture for non-connective algebraic K-theory holds true.

\end{abstract}

\vspace{30pt}

It is known, (see [MV] and [MN)that the Baum-Connes isomorphism conjecture on a given group $G$ follows from the existence of an invertible Dirac morphism associated to that group. In [MN] the authors prove the existence of a Dirac morphism which suffices to construct the assembly map  in the setting of the Baum-Connes conjecture. Here, given a group $G$ we also construct an associated Dirac morphism, and we show that if it is invertible, then the Farrell-jones isomorphism conjecture holds true for that group by using theorem 1 stated below.

\vspace{30pt}

DEFINITION 1:(see [MH], [N] and [AH]) \emph{Given a triangulated category $\cal T$ , we say that $\cal G$ is a set of weak generators of this category if this set is closed under (de-) suspensions and if we get that $X=0$, for $X$ an object of $\cal T$ if for all $G$ in $\cal G$ and all integer $n$ }

\vspace{16pt}

\[ Hom_{\cal T}(\Sigma^{n} G,X)=0 \]

\vspace{16pt}

DEFINITION 2: (see [MH], [SS] and [HPS]) \emph{A stable model category is a pointed stable model category for which the functors $\Sigma$ and $\Omega$ on the homotopy category are invertible.}

\vspace{16pt}

REMARK 1: (see [MH], [SS] and [HPS]) The homotopy category of a stable model category is a triangular category which has all limits and colimits, infinite coproducts and products.

\vspace{16pt}

DEFINITION 2:(see [S]) \emph{Let $\cal T$ be a triangulated category with infinite coproducts. An object C in $\cal T$ is called compact if for any family $(Y_{i})_{i\in I}$ of objects of $\cal T$ one has a natural isomorphism}

\vspace{16pt}

\[\sqcup Hom_{\cal T}(C,Y_{i}) \mapsto Hom_{\cal T}(C,\sqcup Y_{i})\]

\vspace{16pt}

Equivalently any morphism $ C \mapsto \sqcup Y_{i}$ factorizes through a finite subcoproduct.

\vspace {16pt}

DEFINITION 3: (see [N]) \emph{Let $\cal T$ be a triangulated category. A full additive subcategory $\cal S$ is called a triangular subcategory if every object isomorphic to an object of $\cal S$ is in $\cal S$, if $\cal S$ is closed under suspensions and if for any distinguished triangle
}
\vspace{16pt}

\[ X \mapsto Y \mapsto Z \mapsto \Sigma X \]

\vspace{16pt}

\emph{such that the objects $X$ and $Y$ are in $\cal S$, the object $Z$ is also in $\cal S$.
}

\vspace{16pt}

DEFINITION 4: (see [S]) \emph{A triangular subcategory $\cal S$ of a triangular category $\cal T$ is called thick if it is closed under direct summands.}

\vspace{16pt}

DEFINITION 5: (see [S]and [SS])\emph{A triangular subcategory $\cal S$ of a triangular category $\cal T$ is called localizing if it is closed under taking arbitrary coproducts.}

\vspace{16pt}

REMARK 2: By definition, if $\cal S$ is localizing then it is thick.

\vspace{16pt}

DEFINITION 6: (see [SS]) \emph{Given a triangular category $\cal T$ with infinite coproducts, a set $\cal G$ of objects of $\cal T$ is called a generating set of $\cal T$ if the smallest localizing subcategory of $\cal T$ containing $\cal G$ is the whole category $\cal T$.
}
\vspace{16pt}

LEMMA 1: (see [SS] lemma 2.2.1) \emph{Let $\cal T$ be a triangular category with infinite coproducts and let $\cal G$ be a set of compact objects in $\cal T$. Then $\cal G$ is a set of generators of $\cal T$ if and only if $\cal G$ is a set of weak generators of $\cal T$.}

\vspace{16pt}

REMARK 3: (see [R]) By remark 1, the homotopy category $Ho( \cal K)$ of a stable model category $\cal K$ is a triangular category with infinite coproducts, and moreover a set of objects $\cal G$ in this homotopy category is a weak generating set if and only if it is a generating set.

 \vspace{16pt}

 DEFINITION 7: (see [AH] and [N]) i) \emph{An object T in a triangular category $\cal T$ with infinite coproducts is called $\lambda$-small for a regular cardinal $\lambda$, if any map $T \mapsto \sqcup X_{i}$ into an arbitrary coproduct in $\cal T$ factors through some subcoproduct
 }
 \vspace{16pt}

 \[ T \mapsto \sqcup X_{j} \mapsto \sqcup X_{i}\]

 \vspace{16pt}

 \emph{with $J$ a subset of $I$ such that $Card(J)$ is less than  $\lambda$.}

 \vspace{16pt}
 ii) \emph{A set $\cal G$ of objects of $\cal T$ is called $\lambda$-perfect if it satisfies:}

 \vspace{16pt}

 a) $0\in \cal G$

 \vspace{10pt}
 b) \emph{Any map $ G \mapsto \sqcup T_{i}$ with $G$ in $\cal G$, $(T_{i})_{i\in I}$ in $\cal T$ and $card(I)$ less than $\lambda$ factors as}

 \vspace{16pt}

 \[G \mapsto \sqcup G_{i} \mapsto \sqcup f_{i}:\sqcup G_{i} \mapsto \sqcup T_{i}\]

 \vspace{16pt}

 \emph{with $G_{i}$ in $\cal G$ and the maps $f_{i}:G_{i} \mapsto T_{i}$ in $\cal T$.}

 \vspace{16pt}

 iii) \emph{A set $\cal G$ in $\cal T$ is called $\lambda$-compact for a regular cardinal $\lambda$ if every $G$ in $\cal G$ is $\lambda$-small, and also $\cal G$ is $\lambda$-perfect.}

 \vspace{16pt}

 REMARK 4: Any compact object of a triangular category is a $\lambda$-compact object for every regular cardinal $\lambda$.

 \vspace{16pt}

 DEFINITION 8:(see [N],[S] and [AH])  \emph{A triangular category $\cal T$ is called well generated if it has a weak generating set $\cal G$ of $\lambda$-compact objects, for some regular cardinal $\lambda$.}

 \vspace{16pt}

 REMARK 5: In particular, by remark 4, a triangular category $\cal T$ with infinite coproducts and a generating set of compact objects, or equivalently, a weak generating set of compact objects, that is, a compactly generated triangular category, is a well generated category.

\vspace{16pt}

PROPOSITION 1: ([K1]) \emph{A triangular category with infinite coproducts is well generated provided that there is a weak generating set $\cal G$ consisting of $\lambda$-small objects such that for any family of maps $ X_{i} \mapsto Y_{i}$ with $i\in I$ and with induced maps }

\vspace{16pt}

\[Hom_{\cal T} (G,X_{i}) \mapsto Hom_{\cal T}(G,Y_{i})\]

\vspace{16pt}

\emph{being surjective for all $G$ in $\cal G$, the induced map}

\vspace{16pt}

\[ Hom_{\cal T}(G,\sqcup X_{i}) \mapsto Hom_{\cal T}(G,\sqcup Y_{i})\]

\vspace{16pt}

\emph{is also surjective.}

\vspace{16pt}

REMARK 6: The Brown Representability theorem holds for cohomology functors defined over well generated triangular categories (see [N]), but in [K2] it is proven that the Brown Representability theorem holds for cohomology functors under the weaker hypothesis which is the one of proposition 1 but this time deleting the $\lambda$-small condition.

\vspace{16pt}

PROPOSITION 2:(see [AH])\emph{Let $\cal T$ be a well generated triangular category and let $\cal L$ be a localizing subcategory which is generated by a set of objects. Then $\cal L$ is also well generated. Moreover the Verdier quotient $\cal T/\cal L$ is a localization of $\cal T$ and is well generated.
}
\vspace{16pt}

PROPOSITION 3: (see [R]) \emph{Let $\cal K$ be a pointed combinatorial model category, then its homotopy category $Ho(\cal K)$ is well generated.In particular the homotopy category of a stable combinatorial category is well generated.}

\vspace{16pt}

PROPOSITION 4:(see [CGR] Th 3.9)\emph{Let $\cal K$ be a stable combinatorial model category. If Vopenska's principle holds then every localizing subcategory of $Ho(\cal K)$ is single generated.}

\vspace{10pt}

See [AR] for information on the large-cardinal axiom of set theory called Vopenska's principle.

\vspace{16pt}

REMARK 7: It then follows by propositions 2,3 and 4 that under Vopenska's principle, every localizing subcategory of the homotopy category of a stable combinatorial model category is well generated and hence we can apply Brown's representability theorem to that subcategory.

\vspace{30pt}

We are going to apply the definitions, lemmas and propositions stated above to a particular model category, which is $ Spt^{Or(G)}$ See [Hir] or [BM1] for the definition of the model structure on $ Spt^{Or(G)}$ . Here $Spt$ is the model category of spectra of compactly generated Hausdorff spaces with the stable model category structure. $Or(G)$ is the orbit category of a group $G$. The objects are the homogeneous spaces $G/H$ with $H$ a subgroup of $G$ considered as left $G$-sets, and the morphisms are $G$-maps $G/H \mapsto G/K$ given by right multiplication $r_{g}:G/H \mapsto G/K, g_{1}H \mapsto g_{1}gK$ provided $g\in G$ satisfies $g^{-1}Hg \subset K$.

\vspace{2pt}

In particular we need to apply proposition 3 to $Spt^{Or(G)}$, hence we must prove that this category is combinatorial and stable. It is stable since $Spt$ is stable. Recall also that a combinatorial category is a model category which is cofibrantly generated (see [Hir], and is locally presentable, that is, it is cocomplete and accesible, (see [AR]). Since it is a model category, by definition it is cocomplete, so that we must only see that it is accesible and cofibrantly generated. Now $Spt$ is accesible, moreover it is combinatorial (see [R] ex 3.5 iv))and therefore also $Spt^{Or(G)}$ is accesible by ([AR] th2.39, page 96).Also by [BM1 see Th3.5] $Spt^{Or(G)}$ is cofibrantly generated, and hence $Spt^{Or(G)}$ is a combinatorial stable model category and it follows  by proposition 3 that the homotopy category of $Spt^{Or(G)}$ is a well generated triangular category and by remark 7 every localizing subcategory of the homotopy category of $Spt^{Or(G)}$ is also well generated so that we can apply Brown's Representability Theorem to it.

\vspace{30pt}

In [BM2] a Quillen adjunction $ind_{\cal D}^{\cal C}: Spt^{\cal D} \Leftrightarrow  Spt^{\cal C}:res_{\cal D}^{\cal C}$ is defined where $\cal C$ is the orbit category $Or(G)$ and where $\cal D$ stands for the subcategory $Or(G, \cal(VC)$ whose objects $G/H$ are such that $H$ is a virtually cyclic subgroup. There is an induced derived adjunction $Lind_{\cal D}^{\cal C}: Ho(Spt^{\cal D}) \Leftrightarrow Ho( Spt^{\cal C}):res_{\cal D}^{\cal C}$.

\vspace{16pt}

In [LR], page 797, given an associative ring $R$ with unit, a suitable functor $K_{R}:Or(G) \mapsto Spt$ is constructed. This functor has the property that the homotopy groups $\pi_{*}(K_{R}(G/H))$ are canonically isomorphic to the nonconnective algebraic $K$ theory groups $K_{*}(R[H])$ for all subgroup $H$ of $G$, where $R[H]$ is the group ring defined by the group $H$.

\vspace{16pt}

The following theorem is proven in [BM2]

\vspace{10pt}

THEOREM 1: The Farrell-Jones conjecture for nonconnective K-theory is verified for a group $G$ if and only if the image of the functors $K_{R}$ in $Ho( Spt^{\cal C})$ , for all associative rings with unit $R$ belong to $Lind_{\cal D}^{\cal C} (Ho(Spt^{\cal D}))$

\vspace{30pt}

It is known, (see [MV] and [MN)that the Baum-Connes isomorphism conjecture on a given group $G$ follows from the existence of an invertible Dirac morphism associated to that group. In [MN] the authors prove the existence of a Dirac morphism which suffices to construct the assembly map  in the setting of the Baum-Connes conjecture. Here, given a group $G$ we also construct an associated Dirac morphism, and we show that if it is invertible, then the Farrell-jones isomorphism conjecture holds true for that group by using theorem 1.

\vspace{30pt}

DEFINITION 9: \emph{We define by $\cal(CI)$ for the localizing subcategory generated by $Lind_{\cal D}^{\cal C} (Ho(Spt^{\cal D}))$. $\cal(CI)$ is a localizing subcategory of the triangular category $Ho( Spt^{\cal C})$
}

\vspace{25pt}

By analogy with the definition 4.5 in [MN] we define

\vspace{16pt}

DEFINITION 10: Given $X$, an object of $Ho( Spt^{\cal C})$, a $\cal(CI)$ approximation of $X$ is a morphism $f:\widehat X \mapsto X$ with $\widehat X$ belonging to $\cal(CI)$ such that $Lres_{\cal D}^{\cal C}(f)$ is invertible.

\vspace{16pt}

DEFINITION 11 (Dirac morphism): Let $*_{Or(G)}$ be a unit in the symmetric monoidal category $Ho(Spt^{\cal C})$ and use the notation $*_{Or(G)}$ or $*_{\cal C}$ for the sake of simplicity, then a Dirac morphism is a $\cal(CI)$ approximation of $*_{\cal C}$.

\vspace{25pt}

By [CGR] under Vopenska's principle, every localizing subcategory of the homotopy category of a stable combinatorial model category is coreflective. A full subcategory $\cal C$ of a category $\cal T$ is called coreflective if the inclusion $\cal C \hookrightarrow \cal T$ has a right adjoint.The composite $C:\cal T \mapsto \cal T$ is called colocalization onto $\cal C$. $\cal C$ is then the class of $C$-colocal objects. Dually, a reflection $\cal L$ of $\cal T$ is a full subcategory such that the inclusion $\cal L \hookrightarrow \cal T$ has a left adjoint $\cal T \mapsto \cal L$. Then $L: \cal T \mapsto \cal T$ is called localization onto $\cal L$. $\cal L$ is then the class of $L$-local objects. By theorem 1.4 in [CGR] there is a bijective correspondence between coreflections and reflections, and given a coreflection $\cal C$ with colocalization $C: \cal T \mapsto \cal T$ there is a triangle

\vspace{25pt}

\[CX \mapsto X \mapsto LX \mapsto \Sigma(CX)\]

\vspace{25pt}

for all $X$ in $\cal T$, where $L:\cal T \mapsto \cal T$ is the localization associated to the reflection $\cal L$ which corresponds to $\cal C$. Here $\cal L$ coincides with the set of all $X$ in $\cal T$ such that $Hom_{\cal T}(\Sigma^{k}Y,X)=0$ for all $Y$ in $\cal C$ and all integer $k$.

\vspace{25pt}

Observe that our localizing category $\cal(CI)$ is by the results of [CGR] coreflective, therefore for each object $X$ in $Ho(Spt^{Or(G)})$ there is a triangle

\vspace{25pt}

\[P \mapsto X \mapsto N \mapsto \Sigma P\]

\vspace{25pt}

with $P$ belonging to $\cal(CI)$ and $N$ belonging to an orthogonal reflective category (orthogonal in the sense stated above) noted $\cal(CC)$.This an analogous result to the one in [M] Th 70.

\vspace{25pt}

Next we will show that the map $P \mapsto X$ from above is a $\cal (CI)$ approximation of $X$. Note that in that case, by taking $X=*_{\cal C}$ we will have proven the existence of a Dirac morphism.

\vspace{16pt}

We must show by definition that in the triangle above, the map $f:P \mapsto X$ is such that $res_{\cal D}^{\cal C}(f)$ is invertible. Since $res_{\cal D}^{\cal C}$ is a triangulated functor we know that

\vspace{25pt}

\[res_{\cal D}^{\cal C}P \mapsto res_{\cal D}^{\cal C}X \mapsto res_{\cal D}^{\cal C}N \mapsto res_{\cal D}^{\cal C}\Sigma P \]

\vspace{25pt}

is an exact triangle. Hence $res_{\cal D}^{\cal C}(f)$ is invertible if and only if $res_{\cal D}^{\cal C}N=0$ But this fact is immediate since we know that $Hom_{\cal T}(A,N)=0$ for all $A$ in $\cal (CI)$ so that $0=Hom_{\cal T}(A,N)=Hom(Ind_{\cal D}^{\cal C}W,N)=Hom_{\cal T}(W, res_{\cal D}^{\cal C}N)$ implying that $res_{\cal D}^{\cal C}N=0$

\vspace{25pt}

The orthogonality condition $Hom_{\cal T}(A,N)=0$ for all $A$ in $\cal (CI)$ can also be obtained in the following way: Consider the functor $A \mapsto Hom(A,X)$ defined in $\cal (CI)$ and where $X$ in $\cal T$ is fixed. This cohomology functor which takes coproducts into products is defined by what we have already proven, in a well generated localizing category and therefore Brown'srepresentability theorem applies.Hence we have a morphism $f:P \mapsto X$ such that $Hom(A,P) \simeq Hom(A,X)$ where $P$ belongs to $\cal (CI)$. This isomorphism implies by taking $Hom_{\cal T}$ in the triangle

\vspace{25pt}

\[P \mapsto X \mapsto N \mapsto \Sigma P\]

\vspace{25pt}

that $ Hom_{\cal T}(A,N)=0$ for all $A$ in $\cal (CI)$.

\vspace{25pt}

REMARK 8: Observe that if $Y$ belongs to $\cal (CI)$ then taking the product in the symmetric monoidal category $Ho(Spt^{Or(G)})$ with another object $Z$ we get that $Y \otimes Z$ is in $\cal (CI)$. This is immediate if we use the following argument:

\vspace{25pt}

$Hom_{\cal T}(Y\otimes Z,W) \simeq Hom_{\cal T}(Y, Hom_{int}(Z,W))=Hom_{\cal T}(Ind X,Hom_{int}(Z,W))\simeq Hom_{\cal T}(X,Hom_{int}(res Z,res W))\simeq Hom_{\cal T}(X\otimes res Z,res W) \simeq Hom_{\cal T}(Ind (X\otimes res Z),W)$

\vspace{25pt}

where $Hom_{int}$ is the internal $Hom$. Consequently $Y\otimes Z= Ind(X\otimes res Z)$ and we are done.

\vspace{30pt}

By remark 8, note that we can go from the Dirac morphism $D:X \mapsto *_{Or(G)}$ to a $\cal (CI)$ approximation of an object $K_{R}$ with $R$ an associative ring with unit just by taking the product with $K_{R}$ in the Dirac morphism.

\vspace{16pt}

We get $D_{R}:K_{R}\otimes X \mapsto K_{R}$ Therefore if the Dirac morphism $D$ is invertible then so is $D_{R}$ in which case every $K_{R}$ is isomorphic to an object in $\cal (CI)$. If that is so, then by theorem 1, the Farrell-Jones conjecture holds.

\vspace{16pt}

\emph{E-mail address}: valmont8ar@hotmail.com

\end{document}